\documentclass{amsart}
\renewcommand{\i}{\mathrm{i}}
\newcommand{\e}{\mathrm{e}}
\newcommand{\id}{\mathrm{id}}
\renewcommand{\d}{\,\mathrm{d}}
\newcommand{\dd}[1]{\d^2#1}

\newcommand{\R}{\mathbb{R}}
\newcommand{\C}{\mathbb{C}}
\newcommand{\U}{\mathrm{U}}
\newcommand{\Nm}{\mathcal{N}}
\newcommand{\M}{\mathcal M(D)}

\newcommand{\I}{^{-1}}
\newcommand{\0}{^{(0)}}
\newcommand{\Ord}{\mathcal{O}}
\newcommand{\ord}{o}

\DeclareMathOperator{\tr}{tr}
\DeclareMathOperator{\diag}{diag}
\DeclareMathOperator{\RE}{\mathrm{Re}}
\DeclareMathOperator{\IM}{\mathrm{Im}}
\DeclareMathOperator{\im}{im}

\newtheorem{thm}{Theorem}[section]
\newtheorem{lemma}[thm]{Lemma}
\newtheorem{proposition}[thm]{Proposition}
\newtheorem{cor}[thm]{Corollary}
\newenvironment{prf}{\par\noindent{\it Proof:}}{\par\vspace{-\baselineskip}\hfill$\square$\par}
\newenvironment{definition}{\par\noindent{\bf Definition.}}{\par}

\title{Density of eigenvalues of random
normal matrices}
\author{Peter Elbau}
\author{Giovanni Felder}
\address{Department of Mathematics, ETH-Zentrum, 8092 Zurich, Switzerland}
\email{peter.elbau@math.ethz.ch, 
giovanni.felder@math.ethz.ch} 
\begin{document}
\maketitle
\begin{abstract} 
The relation between random normal matrices and
conformal mappings discovered by Wiegmann and Zabrodin is
made rigorous by restricting  normal matrices
to have spectrum in a bounded set. It is shown that for
a suitable class of potentials the asymptotic density
of eigenvalues is uniform with support in the interior
domain of a simple smooth curve.
\end{abstract}

\section{Introduction}
In recent work initiated by P. Wiegmann and A. Zabrodin
 \cite{WZ,KKMWZ,MWZ,KMZ}, a connection between the
normal matrix model
and conformal mappings was discovered. In this model one considers random normal $N\times N$ complex matrices
with probability measure 
\begin{equation}\label{e-pm}
P_N(M)\d M=Z^{-1}_N
\exp\{-(N/t_0)\mathrm{tr}(M^*M-p(M)-p(M)^*)\}\d M,
\end{equation}
where
$\d M$ is a natural measure on the variety of
normal matrices  and $Z_N$ is
the normalization factor.
The result is that, as $N\to\infty$, the density of eigenvalues
is $1/\pi t_0$ times
the characteristic function of a bounded domain in the complex plane.
This domain is characterized by the fact that its exterior harmonic
moments are the coefficients $t_j$ of the polynomial $p$ appearing in
the measure. Moreover, the Riemann mapping of the exterior of the
unit disk onto the exterior of the domain obeys, as a function of
the $t_j$, the equations of the integrable dispersionless Toda hierarchy.

These fascinating results remain however at the level
of formal manipulations of undefined objects, as the
integrals diverge, except in
the simplest case of a polynomial $p$ of degree 2, where the
domain is bounded by an ellipse.

The purpose of this note is to give a setting in which the
above statements make mathematical sense and to give a proof
of these statements.

The problem of divergence of the integral over normal matrices is solved
in a naive way, by restricting the integral to normal matrices whose
eigenvalues lie in a compact domain $D$ of the complex plane.
Then for small $t_0$ the results can be formulated
in terms of {\em polynomial curves}, i.e., curves in the complex plane
admitting a parametrization of the form $w\mapsto h(w)=rw+\sum_{j=0}^n a_jw^{-j}$,
$|w|=1$. For simple polynomial curves the problem of determining the curve out of its
exterior harmonic moments $t_j=\frac1{2\pi\i}\oint_\gamma\bar z z^{-j}\d z$
and the area $\pi t_0$ of the interior domain (the interior domain is
the bounded connected component of the complement of the curve)
has a unique solution for small $t_0$, as we show in Section 5. Our main result
is then:

\begin{thm}\label{t-11} Let $D\subset \mathbb C$ be the closure of a bounded open
set containing the origin.
 Let $p(z)=t_2 z^2+\dots+t_{n+1} z^{n+1}$ be a polynomial 
such that $|z|^2-p(z)-\overline {p(z)} $
has a non-degenerate absolute minimum in $D$ at $z=0$.
 Then there exists a $\delta>0$ so that for all $0<t_0<\delta$,
\begin{enumerate}
\item there exists a unique simple polynomial curve $\gamma$ with exterior
harmonic moments
$t_1=0, t_2,\dots,t_{n+1},0,0,\dots$ and area of interior domain $\pi t_0$;
\item the expectation value of the density of eigenvalues of random normal matrices
with spectrum in $D$ and distribution \eqref{e-pm} converges as $N\to \infty$ to
a uniform distribution with support in the interior domain of $\gamma$.
\end{enumerate}
\end{thm}

The condition on $p$ implies the Hessian condition $|t_2|<\frac12$ and it is fulfilled
if the Hessian condition holds and $t_3,\dots,t_n$ are sufficiently small.
It then follows from results of \cite{TT, WZ} 
 that the curve $\gamma$ as a function of
the $t_j$ in this range provides a solution of the integrable dispersionless 2D Toda hierarchy
obeying the string equation. 

The paper is organized as follows: the basic definitions
of the  random normal matrix model are recalled in Section 2.  We then
introduce the ``equilibrium measure'' as the unique solution of a variational
problem in Section 3 and show in Section 4 that the density of eigenvalues
converges to it. These are either known results or adaptations of
results known for hermitian matrices to our case.
In Section 5 we introduce the notion of polynomial curve and prove Theorem 5.3, which
is a stronger form of part (i) of Theorem \ref{t-11}.
In Section 6 we prove part (ii) of Theorem \ref{t-11}, see Theorem 6.1, and discuss
our results in Section 7.

\section{Eigenvalues of random normal matrices}
We consider the probability measure
\[ P_N(M)\d M=\frac{1}{Z_N}\e^{-N\tr V(M)}\d M,\quad Z_N= \int_{\Nm_N(D)} P_N(M)\d M, \] 
defined by a potential $V$, on the set
\[ \Nm_N(D)=\{M\in\mathrm{Mat}_{\C}(N)\;|\;[M,M^*] = 0,\,\sigma(M)\subset D\} \]
of normal $N\times N$ complex matrices with spectrum in some compact domain $D\subset\C$.

The measure $\d M$ is the Riemannian
volume form on (the smooth part of) $\Nm_N(D)$
with respect to the metric induced from the
standard metric on the vector space $\C^{N^2}$ of all
$N\times N$ matrices. In a parametrization by eigenvalues and unitary matrices it is
given by \cite{LZ}
\[ \d M = \d U\!\!\!\!\!\prod_{1\le i<j\le N}\!\!\!\!\!|z_i-z_j|^2\prod_{i=1}^N\dd{z_i}, \]
where $M=U\diag\big((z_i)_{i=1}^N\big)U^*$ and $\d U$ denotes the normalized $\U(N)$ invariant measure on $\U(N)/\U(1)^N$.

This leads to the probability measure
\begin{align*}
P_N\big((z_i)_{i=1}^N\big)\prod_{i=1}^N\dd z_i &= \frac{1}{Z_N}\e^{-N\sum_{i=1}^N V(z_i)}\!\!\!\!\!\prod_{1\le i<j\le N}\!\!\!\!\!|z_i-z_j|^2\prod_{i=1}^N\dd{z_i}, \\
Z_N &= \int_{D^N} 
\e^{-N\sum_{i=1}^N V(z_i)}\!\!\!\!\!\prod_{1\le i<j\le N}\!\!\!\!\!|z_i-z_j|^2\prod_{i=1}^N\dd{z_i}.
\end{align*}
on the space of eigenvalues $z_i\in D$, $1\le i\le N$.

\section{The equilibrium measure}
We are interested in the behavior of the function $P_N(z)$ as $N\to\infty$. Because the probability that two eigenvalues are equal is always zero, we may consider $P_N$ as function on the set $D_0^N=\{z\in D^N\;|\;z_i\ne z_j\,\forall\,i\ne j\}$.

Introducing the probability measure
\[ \delta_z(A)=\frac{1}{N}\sum_{i=1}^N\chi_A(z_i),\quad z=(z_i)_{i=1}^N, \]
on $D$ ($\chi_A$ shall denote the characteristic function of the set $A$), we write for $z\in D_0^N$
\[ P_N(z) = \frac{1}{Z_N}\exp\left(-N^2\left(\int V(\zeta)\d\delta_z(\zeta)+\iint_{\zeta\ne\xi}\log|\zeta-\xi|\I\d\delta_z(\xi)\d\delta_z(\zeta)\right)\right). \]

Letting $N\to\infty$, only the infimum of the coefficient of $-N^2$,
\begin{align*}
I_0 &:= \!\!\!\inf_{\mu\in\M} I(\mu), \\
I(\mu)&:=\int V(z)\d\mu(z)+\iint_{z\ne\zeta}\log|z-\zeta|\I\d\mu(\zeta)\d\mu(z),
\end{align*}
will be relevant. Here $\mathcal M(D)$ denotes the set of all Borel probability measures on $D$ without point masses (a measure with point masses could only arise from measures $\delta_z$ with some $z_i=z_j$, but then $P_N(z)=0$). Therefore, we safely can neglect the restriction on the double integral.

The precise sense in which the infimum of $I$ controls
the large $N$ behavior of $P_N$ will be discussed in the
next section. Here we consider this reasoning only as a motivation for introducing the variational problem.

Because of $I(\frac{\chi_D}{|D|}\lambda)<\infty$ ($\lambda$ the Lebesgue measure on $\C$ and $|D|=\lambda(D)$ the area of
$D$), $I_0$ is finite.

\medskip
\begin{definition}
An {\em equilibrium measure} for $V$ on $D\subset \C$ is a Borel probability measure $\mu$ on $D$ without point masses so that $I(\mu)=I_0$.
\end{definition}

\begin{thm} Every continuous function $V$ on a 
compact subset $D$ has a unique equilibrium measure.
\end{thm}

In the remaining of this section we prove this theorem, which is a known fact from potential theory, see e.g. \cite{ST},
and give necessary and sufficient conditions for $\mu$
to be an equilibrium measure. The constructions are
adapted from the corresponding results for Hermitian
matrices, see \cite{D}.

\subsection{Existence}
To show the infimum is achieved, we choose a sequence $(\mu_n)_{n=1}^\infty$ in $\M$ with $I(\mu_n)\to I_0$. 
\begin{lemma}\label{l-compactness}
The space of all Borel probability measures on $D$ is sequentially compact.
\end{lemma}
\begin{prf}
By the theorem of Riesz-Markov each Borel measure $\mu$ on $D$ corresponds to exactly one positive, linear functional $\phi_\mu\in C(D)^*$ and by the theorem of Alaoglu, the closed unit-sphere in $C(D)^*$ is weak-*-compact. Therefore, for each sequence $(\mu_n)_{n=1}^\infty\subset\M$, the sequence $(\phi_{\mu_n})_{n=1}^\infty$ contains a weak-*-convergent subsequence $(\phi_{\mu_{n(k)}})_{k=1}^\infty$, i.e.
\[ \exists\,\phi\in C(D)^*\,:\;\phi_{\mu_{n(k)}}(f) \to \phi(f)\quad\forall f\in C(D). \]
Now we find a measure $\mu$ on $D$ with $\phi=\phi_\mu$. This measure fulfills
\[ \int f\d\mu_{n(k)}\to\int f\d\mu\quad(k\to\infty)\quad\forall f\in C(D), \]
and hence is again a Borel probability measure.
\end{prf}
\medskip
Because of this lemma there exists a convergent subsequence $(\mu_{n(k)})$ of $(\mu_n)$ and a Borel probability measure $\mu$ with $\mu_{n(k)}\to\mu$.

To prove that $I(\mu)=I_0$, we estimate with an arbitrary real constant $L$
\begin{align*}
\lim_{k\to\infty}I(\mu_{n(k)}) &= \lim_{k\to\infty}\int V(z)\d\mu_{n(k)}(z)+\lim_{k\to\infty}\iint\log|z-\zeta|\I\d\mu_{n(k)}(\zeta)\d\mu_{n(k)}(z) \\
&\ge\int V(z)\d\mu(z)+\lim_{k\to\infty}\iint\min\{\log|z-\zeta|\I,L\}\d\mu_{n(k)}(\zeta)\d\mu_{n(k)}(z).
\end{align*}
Approximating uniformly the second integrand according to the theorem of Stone-Weierstra\ss\ up to $\varepsilon>0$ with a polynomial in $z$, $\bar z$, $\zeta$ and $\bar\zeta$ and using Fubini's theorem, we get a lower bound for the limit:
\[ \lim_{k\to\infty}I(\mu_{n(k)})\ge\int V(z)\d\mu(z)+\iint\min\{\log|z-\zeta|\I,L\}\d\mu(\zeta)\d\mu(z)-2\varepsilon. \]
Letting first $\varepsilon\to0$ and then $L\to\infty$ shows that $\mu$ has no point masses (otherwise the right hand side would diverge) and $I_0=I(\mu)$.

\subsection{Uniqueness}
Next, we want to show that there is exactly one measure $\mu\in\mathcal M(D)$ with $I(\mu)=I_0$. So suppose $\tilde\mu\in\M$ also fulfills $I(\tilde\mu)=I_0$. Then we consider the family 
\[ \mu_t = t\tilde\mu+(1-t)\mu = \mu+t(\tilde\mu-\mu),\; t\in[0,1], \]
in $\mathcal M(D)$.

\begin{lemma}\label{l-integrability}
Let $\mu$ and $\tilde\mu$ be probability measures such that the function $\log|z-\zeta|\I$ is integrable with respect to $\mu\otimes\mu$ and to $\tilde\mu\otimes\tilde\mu$. Then $\log|z-\zeta|\I$ is also integrable with respect to $\mu\otimes\tilde\mu$.

Additionally, we have the inequality
\begin{equation}\label{e-log}
\iint\log|z-\zeta|\I\d(\tilde\mu-\mu)(\zeta)\d(\tilde\mu-\mu)(z)\ge 0,
\end{equation}
with equality if and only if $\mu=\tilde\mu$.
\end{lemma}
\begin{prf}
We start with the distributional identity
\[ \int\log|z|\I\triangle\varphi(z)\dd z = -2\pi\varphi(0) \]
for any Schwarz function $\varphi$.
Introducing the Fourier transform $\hat\varphi$ of $\varphi$, we find
\begin{align*}
\int\log|z|\I\triangle\varphi(z)\dd z &= -\int\frac{1}{|k|^2}(|k|^2\hat\varphi(k))\dd k \\
&= \frac1{2\pi}\int\frac{1}{|k|^2}\int\triangle\varphi(z)\left(\e^{\frac{\i}{2}(kz+\bar k\bar z)}-f(k)\right)\!\dd z\dd k \\
&= \frac1{2\pi}\iint\frac{1}{|k|^2}\left(\e^{\frac{\i}{2}(kz+\bar k\bar z)}-f(k)\right)\!\dd k\,\triangle\varphi(z)\dd z,
\end{align*}
where $f$ denotes a real, continuous function which is one in the vicinity of zero and becomes zero at infinity. So we have for all $z\in\C\setminus\{0\}$ the equation
\[ \log|z|\I = \frac{1}{2\pi}\int\frac{1}{|k|^2}\left(\e^{\frac{\i}{2}(kz+\bar k\bar z)}-f(k)\right)\!\dd k+C(f) \]
with some real constant $C(f)$.

So, with Tonelli's theorem, we see that
\[ \iint\log|z-\zeta|\I\d(\tilde\mu-\mu)(\zeta)\d(\tilde\mu-\mu)(z)=\frac{1}{2\pi}\int\frac{1}{|k|^2}\left|\int\e^{\frac{\i}{2}(kz+\bar k\bar z)}\d(\tilde\mu-\mu)(z)\right|^2\dd k \]
is non-negative and finite, which immediately implies the integrability of \mbox{$\log|z-w|\I$} with respect to $\mu\otimes\tilde\mu$.
To achieve equality in \eqref{e-log}, we need
\[ \int\e^{\frac{\i}{2}(kz+\bar k\bar z)}\d\mu(z) = \int\e^{\frac{\i}{2}(kz+\bar k\bar z)}\d\tilde\mu(z) \]
for all $k\in\C$, which reads $\mu=\tilde\mu$.
\end{prf}
\medskip

So we can expand $I(\mu_t)$ and obtain
\begin{equation}
\begin{split}\label{e-I expanded}
I(\mu_t) &= I(\mu)+t\int\left(V(z)+2\int\log|z-\zeta|\I\d\mu(\zeta)\right)\d(\tilde\mu-\mu)(z) \\
&\qquad\qquad+t^2\iint\log|z-\zeta|\I\d(\tilde\mu-\mu)(\zeta)\d(\tilde\mu-\mu)(z). 
\end{split}
\end{equation}
Lemma \ref{l-integrability} now states that the coefficient of $t^2$ is non-negative and so the function $t\mapsto I(\mu_t)$ is convex on $[0,1]$. In particular,
\[ I(\mu_t) \le tI(\tilde\mu)+(1-t)I(\mu) = I_0, \]
which implies $I(\mu_t)\equiv I_0$. This requires the last summand in \eqref{e-I expanded} to vanish and so, again with Lemma \ref{l-integrability}, we see that $\mu=\tilde\mu$.

\subsection{A variational form}\label{s-variation}
To determine if a given measure $\mu$ is the equilibrium measure for the potential $V$ on the domain $D$, we may check a variational principle.
\begin{proposition}\label{p-variational principle}
The probability measure $\mu$ is the equilibrium measure for the potential $V$ on the domain $D$ if and only if the function
\begin{equation}\label{e-E}
E(z)=V(z)+2\int\log|z-\zeta|\I\d\mu(\zeta)
\end{equation}
fulfills the relation
\begin{equation}
\int E(z)\d\tilde\mu(z)\ge\int E(z)\d\mu(z)=:E_0\quad\textrm{for all $\tilde\mu\in\mathcal M(D)$.} \label{e-variation 1}
\end{equation}
Additionally, we have the property
\begin{equation}
E(z)\equiv E_0\quad\textrm{$\mu$-almost everywhere.}\label{e-variation 2}
\end{equation}
\end{proposition}

\begin{prf}
Let us first assume $\mu$ is the equilibrium measure. Then, for an arbitrary measure $\tilde\mu\in\mathcal M(D)$, the condition
\[ \frac{\d}{\d t}I(\mu_t)\big|_{t=0} \ge 0,\quad\mu_t=t\tilde\mu+(1-t)\mu,\quad t\in[0,1], \]
has to hold. This means
\[ \int\left(V(z)+2\int\log|z-\zeta|\I\d\mu(\zeta)\right)\d(\tilde\mu-\mu)(z)\ge 0, \]
which immediately implies equation (\ref{e-variation 1}).

Assume on the other side $\mu$ fulfills condition (\ref{e-variation 1}). Then we obtain with the equilibrium measure $\mu_0$
\begin{align*}
I(\mu_0) &= I(\mu)+\int\left(V(z)+2\int\log|z-\zeta|\I\d\mu(\zeta)\right)\d(\mu_0-\mu)(z) \\
&\phantom{=\;I(\mu)}+\iint\log|z-\zeta|\I\d(\mu_0-\mu)(\zeta)\d(\mu_0-\mu)(z) \\
&\ge I(\mu),
\end{align*}
and so $\mu=\mu_0$.

To show the additional statement, we consider for the probability measure $\mu$ the set
\[ B := \{z\in D\;|\;E(z) < E_0\}. \]
If $\mu(B)>0$, the variational principle (\ref{e-variation 1}) for the measure $\tilde\mu = \frac{\chi_B}{\mu(B)}\mu$ would yield
\[ E_0\le\int E(z)\frac{\chi_B(z)}{\mu(B)}\d\mu(z)<E_0, \]
and so $\mu(B)$ has to vanish. Thus, we get condition \eqref{e-variation 2}.
\end{prf}
\medskip

Instead of verifying condition \eqref{e-variation 1}, we prefer to use the following statement:
\begin{cor}\label{c-variational principle}
If for a measure $\mu\in\mathcal M(D)$ the function $E$, defined by equation \eqref{e-E}, fulfills, for some real constant $E_0$, $E(z)\equiv E_0$ on the support of $\mu$ and $E(z)\ge E_0$ everywhere, then $\mu$ is the equilibrium measure.
\end{cor} 

\section{The eigenvalue density}
\begin{definition}
The $k$-point correlation function $R^{(k)}$ is given by 
\[ R^{(k)}_N\big((z_i)_{i=1}^k\big) = \frac{N!}{(N-k)!}\int_{D^{N-k}}P_N\big((z_i)_{i=1}^N\big)\prod_{i=k+1}^N\dd{z_i}. \]
So, the one-point correlation function is up to normalization nothing but the density of the eigenvalues.
\end{definition}
As indicated in the previous section, all the correlation functions can be calculated in the limit $N\to\infty$ out of the equilibrium measure $\mu$, as in the case of hermitian matrix models, see \cite{J,D}.
\begin{thm}\label{t-eigenvalue distribution}
For all $\phi\in C(D^k)$ we have the equality
\begin{equation}\label{e-eigenvalue distribution}
\lim_{N\to\infty}\int_{D^k}\frac1{N^k}\phi\big((z_i)_{i=1}^k\big)R^{(k)}_N\big((z_i)_{i=1}^k\big)\prod_{i=1}^k\dd{z_i}=\int\phi\big((z_i)_{i=1}^k\big)\prod_{i=1}^k\d\mu(z_i).
\end{equation}
I.e. the measure $\frac{1}{N^k}R_N^{(k)}\big((z_i)_{i=1}^k\big)\prod_{i=1}^k\dd{z_i}$ on $D^k$ converges weakly to \mbox{$\prod_{i=1}^k\d\mu(z_i)$.}
\end{thm}
\begin{prf}
Substituting in the left-hand-side of equation \eqref{e-eigenvalue distribution} the definition of the correlation functions and turning our attention on the highest order in $N$, we obtain (because $P_N$ is invariant under the symmetric group)
\begin{align*}
\left<\phi,\frac{1}{N^k}R_N^{(k)}\right>&=\frac{1}{N^k}\int_{D^k}\phi\big((z_i)_{i=1}^k\big)R^{(k)}_N\big((z_i)_{i=1}^k\big)\prod_{i=1}^k\dd{z_i} \\
&= \frac{1}{N^k}\sum_{i_1,\ldots,i_k=1}^N\int_{D^k}\phi\big((z_{i_j})_{j=1}^k\big)P_N\big((z_i)_{i=1}^N\big)\prod_{i=1}^N\dd{z_i}+\ord(1).
\end{align*}

Since for large values of $N$ the probability distribution localizes at values $z\in D_0^N$ with $I(\delta_z)\approx I_0$, let us consider the sets
\[ A_{N,\eta}=\{z\in D_0^N\;|\;I(\delta_z)\le I_0+\eta\},\quad\eta>0. \]
\begin{lemma}
The probability $P_N(D^N\setminus A_{N,\eta})$ drops for $N\to\infty$ exponentially to zero.
\end{lemma}
\begin{prf}
For an absolutely continuous\footnote{This case suffices for our needs, but the restriction is in fact not necessary. Indeed, we could perform an analogous argument for the measures $\d\mu_\varepsilon(z)=\psi_\varepsilon(z)\dd z$, $\psi_\varepsilon(z)=\frac{1}{\pi\varepsilon^2}\int_{B_\varepsilon(z)}\d\mu$. In the limit $\varepsilon\to 0$, where $I(\mu_\varepsilon)\to I_0$, this would yield the desired statement.}
equilibrium measure $\d\mu(z)=\psi(z)\dd z$, we get with Jensen's theorem and $\int I(\delta_z)\prod\d\mu(z_i)=I_0+\ord(1)$
\[ Z_N \ge\int_{\{z\in D^N\,|\,\psi(z_i)\ne 0\;\forall i\}}\hspace{-6em}\e^{-N^2I(\delta_z)-\sum_{i=1}^N\log\psi(z_i)}\prod_{i=1}^N\d\mu(z_i) 
\ge\e^{-N^2I_0+\ord(N^2)}, \]
and therefore,
\[ P_N(D^N\setminus A_{N,\eta})\le\int_{D^N}\e^{N^2I_0+\ord(N^2)-N^2(I_0+\eta)}\prod_{i=1}^N\dd{z_i}=\ord(\e^{-N^2\eta/2}). \]
\end{prf}
\medskip

Let the continuous function $\frac{1}{N^k}\sum\phi\big((z_{i_j})_{j=1}^k\big)$ take its maximum on the compact set $A_{N,\eta}$ at $\zeta$, and set $\nu_{N,\eta} = \delta_\zeta$. Then,
\[ \left<\phi,\frac{1}{N^k}R_N^{(k)}\right>\le\frac{1}{N^k}\sum_{i_1,\ldots,i_k=1}^N\phi\big((\zeta_{i_j})_{j=1}^k\big)=\int\phi\big((z_i)_{i=1}^k\big)\prod_{i=1}^k\d\nu_{N,\eta}(z_i). \]
Because of Lemma \ref{l-compactness}, we find a convergent subsequence $\nu_{N(n),\eta}\to\nu_\eta$ $(n\to\infty)$ with
\[ \varlimsup_{N\to\infty}\left<\phi,\frac{1}{N^k}R_N^{(k)}\right>\le\int\phi\big((z_i)_{i=1}^k\big)\prod_{i=1}^k\d\nu_\eta(z_i). \]

\begin{lemma}
We have $\nu_\eta\in\M$ and, in the limit $\eta\to0$, $I(\nu_\eta)\to I_0$.
\end{lemma}
\begin{prf} 
Using $\zeta\in A_{N(n),\eta}$, we obtain with the cut-off $L\in\R$:
\begin{align*}
I_0+\eta &\ge \frac{1}{N(n)}\sum_{i=1}^{N(n)}V(\zeta_i)+\frac{1}{N(n)^2}\!\!\!\!\sum_{\phantom{(n)}1\le i\ne j\le N(n)}\!\!\!\!\min\{\log|\zeta_i-\zeta_j|\I,L\} \\
&=\int V(z)\d\nu_{N(n),\eta}(z) \\
&\qquad\qquad+\iint\min\{\log|z-\zeta|\I,L\}\d\nu_{N(n),\eta}(\zeta)\d\nu_{N(n),\eta}(z)-\frac{L}{N(n)}.
\end{align*}
Sending first $n$ and then $L$ to infinity brings us to $\nu_\eta\in\M$ and therefore $I_0\le I(\nu_\eta)\le I_0+\eta$.
\end{prf}
\medskip

So, letting $\eta\to0$, a subsequence of $\nu_\eta$ converges to the equilibrium measure $\mu$, and thus,
\[ \varlimsup_{N\to\infty}\left<\phi,\frac{1}{N^k}R_N^{(k)}\right>\le\int\phi\big((z_i)_{i=1}^k\big)\prod_{i=1}^k\d\mu(z_i). \]
Arguing on the same way for the limes inferior concludes the proof.
\end{prf}

\section{Polynomial curves}
\noindent{\bf Definition.} A {\em polynomial curve of degree $n$} is a smooth
simple closed curve
in the complex plane with a parametrization $h:S^1\subset\C\to\C$ of the form
\begin{equation}\label{e-pc1}
h(w)= rw+a_0+a_1w^{-1}+\cdots+a_nw^{-n},\qquad |w|=1,
\end{equation}
with $r>0$ and $a_n\neq0$. The standard (counterclockwise)
orientation of the circle induces an orientation
on the curve. We say that a polynomial curve is positively
oriented if this orientation is counterclockwise, i.e.,
if the tangent vector to the curve makes one full
turn in counterclockwise direction as we go around
the unit circle.

\begin{proposition}\label{p-pc1}
Let $\gamma$ be a positively oriented
polynomial curve with parametrization $h$ of the
form \eqref{e-pc1}. Then $h$, viewed as homolorphic map on $\C^\times$,
restricts to a biholomorphic map from the
exterior of the unit disk onto the exterior of $\gamma$.
\end{proposition}

\begin{prf}
We have to show that $h'(w)\neq 0$ for all $w$ in
the complement of the unit disk. Let $t$ denote the tangent vector
map $w\mapsto t(w)=h'(w)\i w=\i(rw-\sum ja_jw^{-j})$. 
Since $\gamma$ is a simple closed curve,
the map $w\mapsto t(w)/|t(w)|$ is a map of degree 1 from the unit circle
to itself.
 Therefore we have
\begin{eqnarray*}
1&=&\frac1{2\pi}\oint_{|w|=1}\d\,\mathrm{arg}(t(w))\\
&=&\frac1{2\pi\i}\oint_{|w|=1} \frac{t'(w)}{t(w)}\d w\\
&=&N+\frac1{2\pi\i}\oint_{|w|=R} \frac{t'(w)}{t(w)}\d w.
\end{eqnarray*}
Here $N\geq0$ denotes the number of zeros of $t(w)$,
counted with multiplicity, in the
complement of the unit disk and $R$ is so large that it contains them
all. The latter integral is 1 as can be seen by sending $R$ to infinity.
Thus $N=0$, and $h'$ has no zeros in the complement of the unit disk.
\end{prf}
\medskip

A simple consequence of this proposition is that a polynomial curve
is {\em uniquely} parametrized by a map of the form \eqref{e-pc1} with
$r>0$. Indeed, any other conformal mapping of the complement differs
by an automorphism of the complement of the unit disk. But 
non-trivial automorphisms are given by fractional linear transformations
which do not preserve the conditions.



{}From now on, we will only consider polynomial curves encircling
the origin, i.e., such that the origin is contained
in their interior domain. This can always be achieved by a translation,
i.e., a shift of $a_0$.

\medskip
\begin{definition}
The harmonic moments $(t_j)_{j=1}^\infty$ of the exterior domain $D_-$ of a polynomial curve (or more generally of an analytic curve) encircling the origin are defined by
\[ t_j=-\frac{1}{\pi j}\int_{D_-}z^{-j}\dd z=\frac{1}{2\pi\i j}\oint_\gamma\bar zz^{-j}\d z, \]
where only the right integral should be taken as a definition for $j\leq 2$.
\end{definition}

\begin{proposition}\label{p-pc2}
 Let $\gamma$ be a positively oriented
polynomial curve of degree $n$ encircling the origin.
\begin{enumerate}
\item The exterior harmonic moments $t_j$
 of $\gamma$ vanish for all $j>n+1$.
\item
 There exist universal polynomials 
$P_{j,k}\in\mathbb Z[r,a_0,\dots,a_{k-j}]$, 
$1\leq j\leq k$,
so that for $j=1,\dots, n+1$,
\begin{equation}\label{e-harmonic moments}
jt_j=\bar a_{j-1}r^{-j+1}+
\sum_{k=j}^n\bar a_kr^{-k}P_{j,k}(r,a_0,\dots,a_{k-j}).
\end{equation}
Moreover, $P_{j,k}$ is a homogeneous polynomial
of degree $k-j+1$ and
it is also weighted
homogeneous of degree $k-j+1$ for
the assignment $\mathrm{deg}(a_j)=j+1$, $\mathrm{deg}(r)=0$.
\item
The area of the domain enclosed by $\gamma$ is $\pi t_0$ where
\[
t_0=r^2-\sum_{j=1}^nj|a_j|^2.
\]
\end{enumerate}
\end{proposition}

\begin{prf}
Since $\gamma$ encircles the origin,
$h(w)$ never vanishes for $|w|\geq 1$. Hence 
the contour in the formula for $t_j$ may be computed by
taking residues at infinity. For $j\geq1$,
\begin{eqnarray*}
jt_j&=&\frac1{2\pi\i}
\oint_{|w|=1}{\bar h(w^{-1})}{h'(w)}h(w)^{-j}\d w\\
&=&r^{-j}\sum_{k=0}^n\frac{\bar a_k}{2\pi\i}
\oint w^{k-j}
\left(r-\sum_{l=1}^nla_lw^{-l-1}\right)
\left(1+\sum_{l=0}^na_lw^{-l-1}/r\right)^{-j}\d w.
\end{eqnarray*}
The integrals in this sum vanish if 
$k\leq j-2$. The formula for $t_j$ in terms of $a_k$, $r$
is obtained by expanding the geometric
series and picking the coefficient of $w^{-1}$ in the
integrand. This proves (i) and the first part of (ii). The
homogeneity property is clear. The weighted homogeneity
follows by rescaling $w$ in the integral.
The same formula can be used to compute $t_0$, but 
the first term $rw^{-1}$ in $\bar h(w^{-1})$, 
which
does not contribute to the integral and
was omitted for $j\geq 1$, must be added here.
\hfill
\end{prf}
\medskip

\noindent{\it Examples.} The terms in $t_j$ involving
polynomials $P_{j,k}$ with $k\leq3$
are 
\begin{eqnarray*}
{t_{1}}&=&\bar a_{0} - 
{  {r^{-1}
\bar a_{1}{a_{0}}}}  - 
{  {r^{-2\,}\bar a_{2}
(2\,{a_{1}}r - {a_{0}}^{2\,})}}  + { 
 {r^{-3\,}\bar a_{3\,}(3\,{a_{0}}{a_{1}}r - 3\,{a_{2}}r^2
 - {a_{0}}^{3\,})}} +\cdots
\\
{2\,t_{2\,}}&=&{  {r^{-1}\bar a_{1}}}  - 
{  {2\,r^{-2\,}\bar a_{2\,}{a_{0}}}} 
 - {  {3\,r^{-3\,}\bar a_{3\,}({a_{1}}r - {
a_{0}}^{2\,})}} 
 +\cdots
\\
{3\,t_{3\,}}&=&{  {r^{-2\,}\bar a_{2\,}}}  - 
{  {3\,r^{-3\,}\bar a_{3\,}{a_{0}}}}+\cdots 
\end{eqnarray*}

\begin{thm}\label{t-polynomial curve} Let $t_2,\dots,t_{n+1}$ be complex
numbers so that
$|t_2|<1/2$.
Then there exists an $A_0=A_0(t_2,\dots,t_{n+1})>0$ 
so that for all
$A$, $t_1$ with $0<A<A_0$ and $|t_1|^2< A(1/2-|t_2|)$, 
there exists a unique
positively oriented
polynomial curve of degree $\leq n$ encircling the origin,
with area $A$
and exterior harmonic moments $t_1,\dots,t_{n+1},
0,0,\dots$
\end{thm}

\noindent{\it Proof:}
The idea is 
to invert the map $(r,a_0,\dots,a_n)\mapsto (t_0,\dots,t_{n+1})$ 
for small $r$ and $a_0$.
Set $\alpha_j=r^{-j}a_j$, $\rho=r^2$ 
and consider instead the
polynomial map 
\[
F:(\rho,\alpha_0,\dots,\alpha_n)\to (t_0,\dots,t_{n+1}),
\] 
as a map from $\mathbb R\times\C^{n+1}$ to itself.
The first claim is that this map has a smooth 
inverse in some neighborhood of any point 
$t\in\mathbb R\times\C^{n+1}$
such that $t_0=t_1=0$ and $|t_2|\neq1/2$.
By Prop.\ \ref{p-pc2}, this map is given by 
\begin{eqnarray*}
t_0&=&\rho-\sum_{j=1}^n\rho^j j|\alpha_j|^2
\\
jt_j&=&\bar\alpha_{j-1}+\sum_{k=j}^{n}\bar\alpha_k
P_{j,k}(r,\alpha_0,r\alpha_1,\dots,r^{k-j}\alpha_{k-j})
\\
&=&
\bar\alpha_{j-1}+\sum_{k=j}^{n}\bar\alpha_k
P_{j,k}(\rho,\alpha_0,\alpha_1,\dots,\alpha_{k-j}).
\end{eqnarray*}
{}From the contour integral representation of the harmonic moments
we get the integral 
\begin{eqnarray*}
\lefteqn{P_{j,k}(\rho,\alpha_0,\dots,\alpha_{k-j})}\\
&=&
\frac1{2\pi i}\oint_{|w|=R} w^{k-j}\left(1-\sum_{l=1}^nl\alpha_l\rho^lw^{-l-1}\right)
\left(1+\sum_{l=0}^{n}\alpha_l\rho^lw^{-l-1}\right)^{-j}\d w
\end{eqnarray*}
for any sufficiently large $R$. By computing the residue at infinity, we can calculate $P_{j,k}$ and thus $t_j$ 
up to terms of at least second order in
$\alpha_0,\rho$,
\begin{eqnarray*}
t_0&=&\rho\,(1-|\alpha_1|^2)+\cdots,\\
jt_j&=&\bar\alpha_{j-1}-j\,\bar\alpha_j\alpha_0
-(j+1)\rho\,\bar\alpha_{j+1}\alpha_1+\cdots,\qquad j\geq1.
\end{eqnarray*}
Hence $F(0,0,2\,\bar t_2,\dots,(n+1)\,\bar t_{n+1})=
(0,0,t_2,\dots,t_{n+1})$
and the tangent map at this point sends 
$(\dot \rho,\dot\alpha_0,\dots,\dot
\alpha_n)$ to $(\dot t_0,\dots,\dot t_n)$ with
\begin{eqnarray*}
\dot t_0&=&(1-4|t_2|^2)\dot\rho,\\
j\dot t_j&=&\bar{\dot\alpha}_{j-1}-j(j+1)\,
t_{j+1}\dot\alpha_0
-2j(j+1)(j+2)\,{t}_{j+2} t_2\dot\rho,\qquad j\geq1.
\end{eqnarray*}
The tangent map is invertible if $|t_2|\neq 1/2$.
By the inverse function theorem, 
$F$ has a smooth inverse
on some neighborhood of $(0,0,t_2,\dots,t_{n+1})$. If $|t_2|<1/2$, $F$ preserves the positivity
of the first coordinate.

In terms of the original variables, this means that given
any $t=(t_0,t_1,t_2,\dots,t_{n+1})$ with small $t_0>0$,
$t_1$ and such that  $|t_2|\neq 1/2$,  there is a
curve $w\mapsto h(w)$  with $h(w)=rw+\alpha_0+r\alpha_1w^{-1}+\cdots+r^n\alpha_nw^{-n}$ and 
$\alpha_j\simeq(j+1)\,\bar t_{j+1}$. It remains to
show that if $r>0$ is small enough, $h$ parametrizes
a positively oriented
simple closed curve containing the origin
We first show that $h$ is
an immersion. Since $h'(w)=r-r\alpha_1w^{-2}+O(r^2)$
and $\lim_{r\to 0}\alpha_1=2\,\bar t_2$, we see that 
as long as $|t_2|\neq 1/2$, $h'(w)$ does not vanish
on the unit circle. Similarly, we show that $h:S^1\to\C$ 
is injective: we have
\begin{eqnarray*}
|h(w)-h(w')|^2&=&
r|w-w'+2\,\bar t_2(w^{-1}-w'^{-1})|+O(r^2)
\\
&=&
r|w-w'+2\,\bar t_2(\bar w-\bar {w'})|+O(r^2)
\end{eqnarray*}
But the expression in the absolute value
can only vanish for $w\neq w'$ if $|2\,\bar t_2|=1$
which is excluded by the hypothesis. Moreover
$h(w)=rw+\bar t_1+2r\bar t_2w^{-1}+O(r^2)$. 
Therefore $h$ parametrizes a perturbation of
an ellipse centered at $\bar t_1$. The condition
on $t_1$ is a sufficient condition for this
ellipse to contain the origin. 
\hfill $\square$

\medskip

\noindent{\it Example.} Let $k\geq3$ and let us
consider curves with $t_j=0$ for all $j\neq k$. Then
$a_j=0$ for all $j\neq k-1$, so that $h(w)=rw^{-1}+
a_{k-1}w^{k-1}$. The relation between $(t_0,t_k)$ and
$(r,a_{k-1})$ is 
$
t_0=r^2-(k-1)|a_{k-1}|^2$, $k\,t_k=\bar a_{k-1}r^{-k+1}$.
This map is a diffeomorphism from the region
$0<r<(k-1)\,|a_{k-1}|$, which the condition for
 $h(w)$ to be an embedding of the
unit circle,
onto the region 
\[0<t_0<(k(k-1)|t_k|)^{-2/(k-2)}(k-2)/(k-1).\]
As $t_0$ approaches the upper bound for given $t_k$,
the curve develops cusp singularities.

\section{The equilibrium measure for a polynomial curve}
In this section we evaluate the equilibrium measure corresponding to potentials
\begin{equation}\label{e-potential}
V(z) = \frac{1}{t_0}\left(|z|^2-2\RE\sum_{k=1}^{n+1}t_kz^k\right).
\end{equation}
We anticipate the result:
\begin{thm}\label{t-equilibrium measure of a polynomial curve}
For any set $(t_k)_{k=1}^\infty\subset\C$ with $t_1=0$, $|t_2|<\frac12$ and $t_k=0$ for $k>n+1$ and any compact domain $D\subset\C$ containing the origin as interior point and such that $t_0V$, $V$ given by equation \eqref{e-potential}, is positive on $D\setminus\{0\}$, there exists a $\delta>0$ so that for all $0<t_0<\delta$ the equilibrium measure $\mu$ for $V$ on $D$ is given by
\[ \mu=\frac{1}{\pi t_0}\chi_{D_+}\lambda, \]
where $D_+$ denotes the interior domain of the polynomial curve $\gamma$ defined by the harmonic moments $(t_k)_{k=0}^\infty$, and $\lambda$ is 
the Lebesgue measure on $\C$.
\end{thm}
The rest of the section is dedicated to the proof of this theorem.

\subsection{The Schwarz reflection}
We first need the notion of a reflection on an analytic 
 curve (see e.g. \cite{Da}).

\medskip
\begin{definition}
The {\em Schwarz function} of an analytic 
 curve $\gamma$ is defined as the analytic continuation (in a neighbourhood of the curve) of the function $S(z)=\bar z$ on $\gamma$. The {\em Schwarz reflection} $\rho$ for the analytic 
 curve in this domain is the anti-holomorphic map $\rho(z)=\overline{S(z)}$.
\end{definition}
\medskip
\begin{definition}
Under the {\em critical radius} $R$ of the polynomial curve defined by the parametrization $h$ we understand the value
\[ R = \max\{|w|\;|\;h'(w)=0,\; w\in\C\}, \]
which, by definition of the map $h$, is less than 1.
\end{definition}

\begin{lemma}
Let $\gamma$ be a polynomial curve parametrized by $h$ and $R$ its critical radius. Then the Schwarz function $S$ and the Schwarz reflection $\rho$ of the curve $\gamma$ restricted to $h(B_{1/R}\setminus\bar B_R)$, where $B_R$ denotes the open disk with radius $R$ around zero, are biholomorphic respectively anti-biholomorphic maps. They are given by
\begin{equation}
S(z) = \bar h\left(\frac{1}{h\I(z)}\right)\quad\textrm{and}\quad\rho(z)=h\left(\frac{1}{\bar h\I(\bar z)}\right).\label{e-reflection}
\end{equation}
Therefore, $\rho$ maps $\gamma$ identically on itself, $h(B_1\setminus\bar B_R)$ to $h(B_{1/R}\setminus\bar B_1)$ and vice versa. Also, we have $\rho^2=\id$.
\end{lemma}
\begin{prf}
By definition of the Schwarz function, in a neighborhood of $|w|=1$,
\[ S(h(w))=\bar h(w\I). \]
Because for $w\in\C\setminus\bar B_R$ the function $h$ is biholomorphic, we may write
\[ S(z) = \bar h\left(\frac{1}{h\I(z)}\right),\quad \rho(z)=\overline{S(z)}=h\left(\frac{1}{\bar h\I(\bar z)}\right). \]
Taking the derivatives, we find that they do not vanish for $R<|h\I(z)|<\frac1R$.
\end{prf}

\begin{lemma}\label{l-gradient}
In the interior domain $D_+$ of the polynomial curve defined by the parameters $t_k$ from \eqref{e-potential} as its harmonic moments, the function 
\[ E(z)=V(z)+\frac2{\pi t_0}\int_{D_+}\log\left|\frac z\zeta-1\right|\I\dd\zeta \]
is equal to zero and in the exterior domain $D\setminus D_+$, its gradient reads
\begin{equation}\label{e-gradient}
\partial_{\bar z}E(z)=\frac{1}{t_0}(z-\rho(z)).
\end{equation}
\end{lemma}
\begin{prf}
To verify the first statement, we use Green's theorem and obtain
\[ \frac2\pi\int_{D_+}\log|z-\zeta|\I\dd\zeta = -|z|^2+\RE\frac{1}{2\pi\i}\oint_\gamma\left(\log|z-\zeta|\I\bar\zeta+\frac{|\zeta|^2}{\zeta-z}\right)\d\zeta. \]
Integrating by parts of the second integrand yields immediately
\[ \frac2\pi\int_{D_+}\log|z-\zeta|\I\dd\zeta = -|z|^2-2\RE\frac{1}{2\pi\i}\oint_\gamma\log(\zeta-z)\bar\zeta\d\zeta, \]
and expanding the logarithm around $z=0$ leads us to $E(z)=0$ in $D_+$.

For the proof of the second part we write $S=S_\i+S_\e$, where $S_\i$ is analytic 
 in $D_+$ and $S_\e$ in the complement $D\setminus D_+$.
For the exterior function $S_\e$ one finds with the Cauchy integral and Stokes formula:
\[ S_\e(z) = -\frac{1}{2\pi\i}\oint_\gamma\frac{\bar\zeta-S_\i(\zeta)}{\zeta-z}\d\zeta  = \frac{1}{\pi}\int_{D_+}\frac{1}{z-\zeta}\dd\zeta,\quad z\in D\setminus D_+. \]

Because we know $E$ to be constant on $D_+$,
\[ 0 = \partial_zE(z) = \frac{1}{t_0}\left(\bar z-\sum_{k=1}^{n+1}kt_kz^{k-1}-\frac{1}{\pi}\int_{D_+}\frac{1}{z-\zeta}\dd\zeta\right) \]
for all $z\in D_+$. So on $\im\gamma$, and by analytic continuation in the entire domain where $S$ is holomorphic (which includes the exterior domain $D\setminus D_+$), we have $S_\i(z)=\sum_{k=1}^{n+1}kt_kz^{k-1}$.

Therefore, for all $z\in D\setminus D_+$,
\[ \partial_{\bar z}E(z) = \frac{1}{t_0}\left(z-\overline{S_\i(z)}-\overline{S_\e(z)}\right) = \frac{1}{t_0}\left(z-\rho(z)\right). \]
\end{prf}

Let us remark that this proof shows that the Schwarz function $S$ has, at least around infinity, the form
\[ S(z)=\sum_{k=1}^{n+1}kt_kz^{k-1}+\frac{t_0}z+\sum_{k=1}^\infty v_kz^{-k-1}, \]
where the $v_k=\frac1\pi\int_{D_+}z^k\dd z$ denote the harmonic moments of the interior domain $D_+$. This fact was already used in \cite{KKMWZ} to establish a connection between the harmonic moments and the coefficients of the parametrization of $\gamma$.
And as was shown in \cite{KKMWZ}, we find, with the Theorems \ref{t-eigenvalue distribution} and \ref{t-equilibrium measure of a polynomial curve}, that the $v_k$ are nothing but the expectation values of $\frac{t_0}N\tr(M^k)$ with respect to the probability measure $P_N(M)\d M$ in the limit $N\to\infty$ and hence, are completely determined by the harmonic moments of the exterior domain and $t_0$.

\subsection{The Gaussian case}
In this case, where the polynomial curve is an ellipse, we are able to calculate the equilibrium measure on $\C$ explicitly. 
\begin{proposition}
The equilibrium measure $\mu$ for the potential
\[ V(z)=\frac{1}{t_0}(|z|^2-t_2z^2-\bar t_2\bar z^2),\quad |t_2|<\frac12, \]
on $\C$ is
\[ \mu = \frac{1}{ab}\chi_{D_+}\lambda, \]
where $D_+$ denotes the interior of the ellipse
\begin{equation}\label{e-ellipse}
\frac{\RE(\sqrt{t_2}z)^2}{a^2}+\frac{\IM(\sqrt{t_2}z)^2}{b^2}=|t_2|,\quad a=\sqrt{\frac{1+2|t_2|}{1-2|t_2|}t_0},\quad b=\sqrt{\frac{1-2|t_2|}{1+2|t_2|}t_0}.
\end{equation}
\end{proposition}
\begin{prf}
As polynomial curve, the ellipse (\ref{e-ellipse}) has the parametrization
\begin{equation}\label{e-param. ell.}
h(w) = r(w+2t_2w\I),\quad r=\frac{a+b}2.
\end{equation}
We check that the given measure $\mu$ is the equilibrium measure by verifying the conditions of Corollary \ref{c-variational principle}. To this end, let us introduce for $|w|>1$ the function
\[ \mathcal E(w) = E(h(w)),\quad E(z) = V(z)+\frac2{\pi t_0}\int_{D_+}\log\left|\frac z{\zeta}-1\right|\I\dd\zeta. \]
Integrating equation (\ref{e-gradient}) and its complex conjugate analog we get for $\mathcal E(w)$ the expression 
\[ \mathcal E(w)=\frac1{t_0}\left(|h(w)|^2-|h(1)|^2-2\RE\int_1^w\bar h(w\I)h'(w)\d w\right). \]
Substituting in this expression relation (\ref{e-param. ell.}) for $h$, we obtain
\begin{align*}
t_0\mathcal E(w) &= r^2\big(|w+2t_2w\I|^2-|1+2t_2|^2\big) \\
&\qquad\qquad-2r^2\RE(\bar t_2w^2+(1-4|t_2|^2)\log w+t_2w^{-2}-(t_2+\bar t_2)) \\
&= r^2\big(|w|^2-1-4|t_2|^2+4|t_2|^2|w|^{-2}+(1-4|t_2|^2)\log(|w|^{-2})\big) \\
&\qquad\qquad+2r^2\RE(2t_2\bar ww\I-t_2(\bar ww\I)|w|^{-2}-t_2(\bar ww\I)|w|^2) \\
&=r^2(|w|^2-1)(1-4|t_2|^2|w|^{-2})+r^2(1-4|t_2|^2)\log(|w|^{-2}) \\
&\qquad\qquad-2r^2(|w|^2-1)(1-|w|^{-2})\RE(t_2\bar ww\I) \\
&=r^2(|w|^2-1)(1-2|t_2|)(1+2|t_2||w|^{-2})+r^2(1-4|t_2|^2)\log(|w|^{-2}) \\
&\qquad\qquad+2r^2(|w|^2-1)(1-|w|^{-2})(|t_2|-\RE(t_2\bar ww\I)).
\end{align*}
We are now ready to start estimating $\mathcal E(w)$ for $|w|>1$ and $|t_2|<\frac12$. The last bracket we can estimate by
\[ |t_2|-\RE(t_2\bar ww\I)\ge |t_2|-|t_2\bar ww\I|=0. \]
It therefore remains to show
\[ (|w|^2-1)(1+2|t_2||w|^{-2})+(1+2|t_2|^2)\log(|w|^{-2})\ge 0, \]
which follows immediately out of the following lemma.
\begin{lemma}
For $0\le\alpha\le 1$ the function 
\[ f(x) = (x-1)(1+\alpha x\I)-(1+\alpha)\log x \]
is non-negative on the interval $[1,\infty)$.
\end{lemma}
\begin{prf}
For the function $f$ and its derivatives $f'$ and $f''$ we have
\[ f(1)=0,\quad f'(1)=0\quad\textrm{and}\quad f''(x)=\frac{1}{x^2}\left(1+\alpha-\frac{2\alpha}{x}\right)\ge0. \]
\end{prf}
\medskip

So we showed that $E(z)\ge0$ for all $z\in\C\setminus D_+$. Because of Lemma \ref{l-gradient}, we also know that in the interior domain $D_+$, $E$ is zero and we therefore can apply Corollary \ref{c-variational principle} to see that $\mu$ is indeed the equilibrium measure for the potential $V$ on $\C$. 
\end{prf}

\subsection{The proof of theorem \ref{t-equilibrium measure of a polynomial curve}} 
As in the Gaussian case we are going to show that the measure $\mu$ given in the theorem fulfills the conditions of Corollary \ref{c-variational principle} and is therefore the uniquely defined equilibrium measure.

Because Theorem \ref{t-equilibrium measure of a polynomial curve} is only valid for interior domains with small area, we are going to consider the asymptotical behaviour $t_0\to0$, where the harmonic moments $(t_k)_{k=1}^\infty$ are kept fixed. 
To catch the asymptotical behavior of the corresponding polynomial curve, let us parametrize it as in the proof of Theorem \ref{t-polynomial curve}:
\[ h(w) = rw+\sum_{j=0}^nr^j\alpha_jw^{-j}. \]
Then, for $r\to0$, we have $r^2\simeq t_0$, $\alpha_j\simeq(j+1)\bar t_{j+1}$, $j\ge1$, and, because we set $t_1=0$, $\alpha_0\simeq r^2$.

\begin{lemma}
The critical radius of $h$ is asymptotically constant for $r\to 0$:
\[ R=\sqrt{|\alpha_1|}+\Ord(r). \]
\end{lemma}
\begin{prf}
The roots of the function $h'(w)=r-r\alpha_1w^{-2}-\ldots-nr^n\alpha_nw^{-n-1}$ are in zeroth order at $\pm\sqrt{\alpha_1}$ and ($n-1$-times degenerated) at zero.
\end{prf}
\medskip

We consider now for $z\in D$ respectively for $w\in h\I(D\setminus D_+)$ the functions
\begin{align}
E(z) &= V(z)+\int_{D_+}\log\left|\frac z\zeta-1\right|\I\dd\zeta\quad\textrm{and} \label{e:energy} \\
\mathcal E(w) &= E(h(w)) = \frac1{t_0}\left(|h(w)|^2-|h(1)|^2+2\RE\int_1^w\bar h(\tilde w\I)h'(\tilde w)\d\tilde w\right).
\end{align}

We already showed in Lemma \ref{l-gradient} that $E(z)\equiv 0$ in $D_+$, so the first condition of the corollary is satisfied (this is essentially the way we have chosen our potential $V$).

Now, also with Lemma \ref{l-gradient}, we see that $E\ge0$ in the vicinity of the curve $\gamma$, strictly speaking in the domain $h(B_{1/R}\setminus\bar B_1)$. Indeed, if we look at the connected components of the contour lines of the function $E$ (which are smooth curves in the considered domain because there $\partial_{\bar z}E\ne 0$), we see that the gradient vector $\partial_{\bar z}E$ always points outwards, i.e. in the exterior domain of the contour line. Therefore, the value of $E$ on the contour lines is increasing outwards as desired.

A bit farther from the curve, i.e. for $1/R\le|w|<r^{-\alpha}$, $0<\alpha<\frac13$, the function $\mathcal E$ equals asymptotically the one of the corresponding ellipse $h\0(w) = rw+\alpha_0+r\alpha_1w\I$, we denote it by $\mathcal E\0$. Indeed, remarking that for the area of this ellipse we have $t_0\0 = t_0+\Ord(r^4)$ and that $\mathcal E\0(w)=\Ord(r^{-2\alpha})$, we obtain, uniformly in $w$,
\begin{align*}
\mathcal E(w)-\mathcal E\0(w) &= \frac1{t_0}(|h(w)|^2-|h\0(w)|^2-|h(1)|^2+|h\0(1)|^2) \\
&\qquad+\frac2{t_0}\RE\int_1^w(h(\tilde w\I)-h\0(\tilde w\I))h'(\tilde w)\d\tilde w \\
&\qquad+\frac2{t_0}\RE\int_1^w(h'(\tilde w)-{h\0}'(\tilde w))h\0(\tilde w\I)\d\tilde w \\
&\qquad+\frac1{t_0}\left(t_0\0-t_0\right)\mathcal E\0(w) \\
&= \Ord(r^{1-\alpha})+\Ord(r^{1-3\alpha})+\Ord(r^{1-2\alpha})+\Ord(r^{2-2\alpha}) \\
&\to0\quad(r\to0)
\end{align*}
Because $\mathcal E\0(w)\ge C>0$ for all $|w|\ge \frac1R$ and $r>0$, we may choose $r$ so small that $|\mathcal E(w)-\mathcal E\0(w)|<\mathcal E\0(w)$ and so $\mathcal E(w)>0$ for all $w\in B_{r^{-\alpha}}\setminus B_{1/R}$.

It remains the domain, where $|w|\ge r^{-\alpha}$. 
For $k\ge2$, we obtain
\begin{align*}
h(w)^k-(rw)^k &= \sum_{l=1}^k\begin{pmatrix}k\\l\end{pmatrix}(rw)^{k-l}\left(\sum_{j=0}^nr^j\alpha_jw^{-j}\right)^l \\
&=\sum_{l=1}^k\begin{pmatrix}k\\l\end{pmatrix}r^kw^{k-2l}\left(\sum_{j=0}^nr^{j-1}\alpha_jw^{1-j}\right)^l 
=\Ord(r^2).
\end{align*}
And since 
\begin{align*}
t_0V(z)&=|z|^2-t_2z^2-\bar t_2\bar z^2+\ord(|z|^2) \\
&=\frac12\left(|z-2t_2\bar z|^2+(1-4|t_2|^2)|z|^2\right)+\ord(|z|^2)
\end{align*}
for $z\to 0$ and $V>0$ on $D\setminus\{0\}$, we may find a constant $C>0$ such that $t_0V(z)\ge C|z|^2$ for all $z$ in the compact domain $D$.

Therefore, for $r\to0$,
\[ V(h(w))=V(rw)+\Ord(1)\ge\frac C{t_0}|rw|^2+\Ord(1), \]
which tends to infinity at least as $r^{-2\alpha}$.

On the other side, the integral over the logarithm in \eqref{e:energy} diverges for $r\to0$ only as $\log r$. So we have for $r$ small enough $E(h(w))>0$ for all $w\in h\I(D\setminus D_+)\setminus B_{r^{-\alpha}}$.

This proves now $E(z)\ge0$ for all $z\in D$ and therefore, with Corollary \ref{c-variational principle}, that $\mu$ is the equilibrium measure.

\subsection{Shifting the origin}
As a little generalization, we consider the case where $t_1\ne0$. This corresponds to a shift of the origin. Therefore, we like to define the harmonic moments also for a curve which does not encircle the origin.

\medskip
\begin{definition}
Let
\[ h(w)=rw+\sum_{j=0}^na_jw^{-j} \]
parametrize a polynomial curve of degree $n$. Then the harmonic moments $(t_k)_{k=1}^{n+1}$ are given by the equation system \eqref{e-harmonic moments}. All other harmonic moments are set to zero.
\end{definition}
Proposition \ref{p-pc2} tells us that this definition coincides with the previous one if the origin is in the interior domain of the curve.
\begin{cor}
For any set $(t_k)_{k=1}^\infty\subset\C$ with $|t_2|<\frac12$ and $t_k=0$ for $k>n+1$ and any compact domain $D\subset\C$ such that $U(z)=|z|^2-2\RE\sum_{k=1}^{n+1}t_kz^k$ has exactly one absolute minimum in the interior of $D$, there exists a $\delta>0$ so that for all $0<t_0<\delta$ the equilibrium measure $\mu$ for $V=\frac1{t_0}U$ on $D$ is given by
\[ \mu=\frac{1}{\pi t_0}\chi_{D_+}\lambda, \]
where $D_+$ denotes the interior domain of the polynomial curve $\gamma$ defined by the harmonic moments $(t_k)_{k=0}^\infty$.
\end{cor}
\begin{prf}
Let us first shift the origin by $a_0$, such that $V(z+a_0)$ has its absolute minimum in $0$. Thereby, the potential gets the form
\[ V(z+a_0)=\frac1{t_0}\left(|z|^2-2\RE\sum_{k=2}^{n+1}t'_kz^k\right)+V(a_0), \]
where the $t'_k$ are the harmonic moments of the shifted curve $\gamma-a_0$. Indeed, the coefficients of $V$ and the harmonic moments depend polynomially on the shift $a_0$. Because we know them to coincide as long as the origin is inside $D_+$, they do so for all $a_0$. 

Applying now Theorem \ref{t-equilibrium measure of a polynomial curve} for the shifted potential gives the desired result.
\end{prf}

\section{Discussion}
We have shown that under suitable assumptions on the polynomial $p$ appearing
in the potential and on the integration range $D$ of the eigenvalues,
the asymptotic density is uniform
with support on a domain uniquely determined by the coefficients of the polynomial $p$.

It would be interesting to understand what happens at the range of validity of our
assumptions. If the potential $V$ has more than one minimum in $D$ then one should
expect for small $t_0$ to have an equilibrium measure with disconnected support, so
that a description by a polynomial curve cannot be valid. Also as $t_0$ gets bigger,
polynomial curves develop singularities and become non-simple. The question is then
what happens to the eigenvalues. Finally we note that polynomial curves are (real sections of
complex) rational curves.
Curves of higher genus should arise by replacing $p$ by more general holomorphic
functions.

\end{document}